\documentstyle[pb-diagram,leqno,12pt]{article}
\input epsf
\input amssym

\def\di{\displaystyle}
\def\be{\begin{eqnarray}}
\def\ee{\end{eqnarray}}

\begin{document}
\date{}
\title{\sc Induced Invariant Finsler Metrics on Quotient Groups }
\author{ E. Esrafilian\\ and\\ \vspace{1cm} H. R. Salimi Moghaddam\\
\emph{Department of Mathematics},\\ \emph{University of Isfahan}, \\ \emph{Isfahan,81746-73441-Iran.} \\
e-mail:\\
 salimi.moghaddam@gmail.com and hr.salimi@sci.ui.ac.ir}

\maketitle \abstract{} \setcounter{equation}{0} In this paper we
show that every invariant Finsler metric on Lie group $G$,
induces an invariant Finsler metric on quotient group $G/H$ in
the natural way, where $H$ is a closed normal Lie subgroup of $G$.



\section{Introduction.}
\setcounter{equation}{0} The study of invariant structures on
homogeneous manifolds is an important problem in geometry. K.
Nomizu obtained many interesting properties of invariant
Riemannian metrics on homogeneous space $G/H$. He introduced
\textit{reductive} homogeneous spaces and studied invariant
Riemannian metrics and the existence and properties of invariant
affine connections on reductive homogeneous spaces (See \cite{[4]}
and \cite{[6]}). Also some curvature properties of invariant
Riemannian metrics on Lie groups and homogeneous spaces have
studied by J. Milnor and H. Samelson (See \cite{[5]} and
\cite{[7]}). So it is important to study invariant Finsler metrics
which are a generalization of invariant
Riemannian metrics.\\
Some properties of invariant Finsler metrics on reductive
homogeneous manifolds are studied in \cite{[2]} and \cite{[3]} by
S. Deng and Z. Hou. The authors of these papers obtained a
necessary and sufficient condition for homogeneous manifolds to
have invariant Finsler metrics. Then they studied bi-invariant
Finsler metrics on Lie groups and obtained a necessary and
sufficient condition for a Lie group to have bi-invariant Finsler
metrics.\\
In this paper we show that every invariant Finsler metric on a Lie
group $G$ induces an invariant Finsler metric on quotient group
$G/H$ in the natural way, where $H$ is a closed normal Lie
subgroup of $G$.\\

\textbf{Note.} In this article we \textbf{do not assume} the
quotient groups $G/H$ are reductive.

\section{Preliminaries.}
\setcounter{equation}{0}
\paragraph{2.1. Definition.} A Minkowski norm on ${\Bbb{R}}^n$ is
a nonnegative function ${\Bbb{F}}:{\Bbb{R}}^n\rightarrow
[0,\infty)$ which has the following properties:\\

(i) ${\Bbb{F}}$ is ${\mathcal{C}}^\infty$ on ${\Bbb{R}}^n \setminus \{0\}$.\\

(ii) ${\Bbb{F}}(\lambda y)=\lambda {\Bbb{F}}(y)$ for all $\lambda > 0$ and $y\in {\Bbb{R}}^n$\\

(iii) The $n\times n$ matrix $(g_{ij})$, where $g_{ij}(y):=
[\di\frac{1}{2}{\Bbb{F}}^2]_{y^iy^j}(y)$, is
positive-\hspace*{1cm}definite at all $y\neq 0$.
\paragraph{2.2. Definition.} Let $M$ be an $n-$dimensional smooth manifold. Also let $TM$ be
the tangent bundle of $M$. A function $F:TM\rightarrow
[0,\infty)$ is called a Finsler metric if it has the following properties:\\

(i) $F$ is ${\mathcal{C}}^\infty$ on the slit tangent bundle $TM\setminus 0$.\\

(ii) For each $x\in M$, ${\Bbb{F}}_x:=F|_{T_xM}$ is a Minkowski
norm on $T_xM$.\\

If the Minkowski norm satisfies ${\Bbb{F}}(-y)= {\Bbb{F}}(y)$,
then one has the absolutely homogeneity $F(\lambda y)= |\lambda|
F(y)$, for any $\lambda\in{\Bbb{R}}$. Every absolutely
homogeneous Minkowski norm is a norm in the sense of functional
analysis.\\
Every Riemannian manifold $(M,g)$ by defining \be
F(x,y):=\sqrt{g_x(y,y)} \hspace*{2cm}x\in M, y\in T_xM\nonumber\ee
is a Finsler manifold (For more details about Finsler geometry see
\cite{[1]}).

We also use the following notations:\\
\begin{itemize}
  \item $R_g:G\rightarrow G$, right translation, $R_g(h)=hg$.
  \item $L_g:G\rightarrow G$, left translation, $L_g(h)=gh$.
  \item $\nu:G\rightarrow G$, inversion, $\nu(g)=g^{-1}$.
  \item $e\in G$, the unit element.
\end{itemize}
We use $F$ for Finsler metrics on Lie group $G$, ${\Bbb{F}}$ for
Minkowski norms on a specific tangent space $T_xM$ or a real
vector space ${\Bbb{R}}^n$ and ${\mathcal{F}}$ for Finsler metrics
on quotient group $G/H$. Also if $f:M\rightarrow N$ is a smooth
function between manifolds and $x\in M$, we denote by
$T_xf:T_xM\rightarrow T_{f(x)}N$ the derivative of $f$ at $x$. If
$f:M\rightarrow N$ is a local diffeomorphism then $T_xf$ is an
isomorphism of vector spaces, yielding for each vector field
$Y\in{\mathcal{X}}(N)$ on $N$ a vector field
$f^{\ast}Y\in{\mathcal{X}}(M)$ defined by
$(f^{\ast}Y)(x)=(T_xf)^{-1}Y(f(x))$.

\section{Induced Invariant Finsler Metrics on \\ {\hspace*{2.5cm} Quotient Groups.}}
\setcounter{equation}{0} Let $G$ be a compact connected Lie group,
$H$ a closed subgroup of $G$, $M=G/H$ the homogeneous space which
consists of left cosets of $zH$, $z\in G$, and $p:G\rightarrow M$
be the natural projection of $G$ onto $M$. The group $G$ admits a
bi-invariant Riemannian metric. Now we can to obtain a Riemannian
metric on $M$ in the following way which is invariant under the
customary action of $G$ on $M$:\\

Let $x\in M$, $X\in T_xM$, $z\in G$ and $Z\in T_zG$, such that
$p(z)=x$, $T_z p(Z)=X$ and let $Z$ be orthogonal to the coset
$z.H$ (a submanifold of G) at $z$. Now we define
$|X|:=|Z|$ (see\cite{[7]}).\\

But in Finsler geometry we have no orthogonality for tangent
vectors so we can't to use the above way. In this article we try
to replace the orthogonality condition, by other conditions such
that by define ${\mathcal{F}}(X):=F(Z)$, have a bi-invariant
Finsler metric on $M$.\\
From now $G$ is an arbitrary finite dimensional Lie group (no
necessarily compact or connected).

For construct a left or right invariant Finsler metric on a Lie
group $G$, it is sufficient to have a Minkowski norm on $T_eG$
such as ${\Bbb{F}}_0$, then
define \be \hspace*{3cm}&F:TG\rightarrow [0,\infty)& \nonumber \\
&F(x,y)={\Bbb{F}}_0(T_xL_{x^{-1}}y)& \hspace*{1.5cm} x\in G, y\in
T_xG \nonumber \ee for left invariant Finsler metrics, and \be
\hspace*{-0.5cm}F(x,y)={\Bbb{F}}_0(T_xR_{x^{-1}}y) \nonumber \ee
for right invariant Finsler metrics.
\paragraph{3.1. Lemma.}
Assume that $G$ is any Lie group and $H$ any closed subgroup, and
denote by ${\frak{g}}$ and ${\frak{h}}$ the Lie algebras of right
invariant vector fields of $G$ and $H$, respectively. Let $V$ be
a vector subspace complementary to ${\frak{h}}$ in ${\frak{g}}$,
that is, ${\frak{g}}=V\bigoplus {\frak{h}}$, and $M:=G/H$ be the
quotient manifold consists of left cosets $zH$, $z\in G$. Then
$\pi: V\rightarrow {\mathcal{X}}(M)$ defined unambiguously by
$\pi(X)(p(z))=T_zp(X(z))$ is a linear function, where
$p:G\rightarrow M:=G/H$ is the natural projection.\\

\noindent\textit{Proof}: Since $\pi$ defined by $Tp$ so $\pi$ is
linear. Let $\{X_1,\cdots,X_k,X_{k+1},\cdots,\\ X_n\}$ be a basis
of Lie algebra of the Lie group $G$ (consists of right invariant
vector fields) such that $\{X_1,\cdots,X_k\}$ is a basis of the
Lie algebra of closed Lie subgroup $H$. So
$\{X_{k+1},\cdots,X_n\}$ is a basis of vector space $V$. We must
show $\pi$ is
welldefined.\\
Assume that $z_1, z_2\in G$ and $p(z_1)=p(z_2)=x$, therefore
$z_1H=z_2H$, so $z_1^{-1}z_2\in H$.\\
Also for $h\in H$ we have $p\circ R_h=p$ because for any $g\in G $
\be p\circ R_h(g)=p(gh)=ghH=gH=p(g).\nonumber\ee So $p\circ
R_{z_1^{-1}z_2}=p$.\\
Let $f \in {\mathcal{C}}^{\infty}(M,\Bbb{R})$ be a real valued
differentiable function, then by attention to the fact that (for
$i=1,\cdots,n$)$X_i$ is right invariant we have \be
(\pi(X_i)(p(z_2)))f&=&(T_{z_2}p(X_i(z_2)))f\nonumber\\
            &=&(T_{z_2}p(T_{z_1}R_{z_1^{-1}z_2}(X_i(z_1))))f\nonumber\\
            &=&X_i(z_1)(f\circ p\circ R_{z_1^{-1}z_2})\nonumber\\
            &=&X_i(z_1)(f\circ p)\nonumber\\
            &=&(T_{z_1}p(X_i(z_1)))f\nonumber\\
            &=&(\pi(X_i)(p(z_1)))f.\nonumber
            \ee
So for any $X\in V$ such that $X=\sum_{i=k+1}^n\lambda_iX_i$ we
have \be (\pi(X)(p(z_2)))f=(\pi(X)(p(z_1)))f\nonumber\ee
Therefore the definition of $\pi$ is welldefined.\hfill\ $\Box$
\paragraph{3.2. Lemma.} Consider the assumptions of Lemma 3.1 and
also suppose that $H$ is a closed normal Lie subgroup of $G$.
Then \be T_zp:V(z)\rightarrow T_{p(z)}M \nonumber\ee is an
isomorphism, where $V(z)=$span$\{X_{k+1}(z),\cdots,X_n(z)\}$.\\

\noindent\textit{Proof}: For $i=1,\cdots, k$ we have $X_i(e)\in
T_eH$, and also \be R_z:H\rightarrow Hz=zH \nonumber\ee is a
diffeomorphism, so $X_i(z)=T_eR_zX_i(e)\in T_zzH$. Therefore\\
$X_i(z)\in\ker(T_zp:T_zG\rightarrow T_{p(z)}M)$. But we know that
$T_zp:T_zG/T_zzH\rightarrow T_{p(z)}M$ is an isomorphism of
vector spaces and $T_zG/T_zzH\simeq V(z)$, so \be
T_zp:V(z)\rightarrow T_{p(z)}M \nonumber\ee is an isomorphism of
vector spaces.\hfill\ $\Box$
\paragraph{3.3. Theorem.} Assume that $G$ is any $n-$dimensional
Lie group, $H$ any closed normal Lie subgroup, $M=G/H$ the
quotient group and $p:G\rightarrow M$ is the natural projection.
If $F$ is a right invariant Finsler metric on $G$, then there is
a Finsler metric on $M$ induced by $F$ such that is invariant
under the natural right action of $G$ on $M$.\\

\noindent \textit{Proof}: Suppose that $\frak{g}$ and $\frak{h}$
are the algebras of right invariant vector fields of $G$ and $H$,
respectively, and $\{X_1,\cdots,X_k,X_{k+1},\cdots,X_n\}$ is a
basis of $\frak{g}$ such that $\{X_1,\cdots,X_k\}$ is a basis of
$\frak{h}$. Let $V$ be a vector subspace complementary to
$\frak{h}$ in $\frak{g}$, that is, $\frak{g}=V\bigoplus \frak{h}$.
Assume that $x\in M$ and $X\in T_xM$ is a tangent vector at $x$.
Let $z\in G$ and $Z\in V(z)$ such that $p(z)=x$ and $T_zp(Z)=X$.
(By Lemma 3.2 for any fixed $z$ such that $p(z)=x$, there is a
unique $Z\in V(z)$ such that $T_zp(Z)=X$) In this situation we
define \be {\mathcal{F}}(X):=F(Z) \nonumber\ee At the first we
show that this definition is welldefined.

For this, we must to show that the definition of ${\mathcal{F}}$
is independent of choice of $z$. \\
Assume $z_1,z_2\in G, Z_1\in V(z_1), Z_2\in V(z_2)$ and
$p(z_1)=p(z_2)=x,\\ T_{z_1}p(Z_1)=T_{z_2}p(Z_2)=X,
Z_1=\sum_{i=k+1}^n \lambda_i X_i(z_1), Z_2=\sum_{i=k+1}^n \mu_i
X_i(z_2)$.\\
Now we can write \be
T_{z_1}p(T_{z_2}R_{z_2^{-1}z_1}(Z_2))&=&T_{z_1}p(T_{z_2}R_{z_2^{-1}z_1}(\sum_{i=k+1}^n\mu_i
                                        X_i (z_2)))\nonumber\\
                                     &=&T_{z_1}p(\sum_{i=k+1}^n\mu_iT_{z_2}R_{z_2^{-1}z_1}(X_i(z_2)))\nonumber\\
                                     &=&T_{z_1}p(\sum_{i=k+1}^n\mu_i(X_i(z_1))\nonumber\\
                                     &=&\sum_{i=k+1}^n\mu_iT_{z_1}p(X_i(z_1))\nonumber\\
                                     &=&\sum_{i=k+1}^n\mu_i(\pi(X_i)(p(z_1)))\nonumber\\
                                     &=&\sum_{i=k+1}^n\mu_i(\pi(X_i)(p(z_2)))\nonumber\\
                                     &=&\sum_{i=k+1}^n\mu_iT_{z_2}p(X_i(z_2))\nonumber\\
                                     &=&T_{z_2}p(\sum_{i=k+1}^n\mu_i(X_i(z_2)))\nonumber\\
                                     &=&T_{z_2}p(Z_2)=X\nonumber\ee
But since $T_{z_1}p:V(z_1)\rightarrow T_xM$ is an isomorphism of
vector spaces, we have $T_{z_2}R_{z_2^{-1}z_1}(Z_2)=Z_1$ (This
also shows that for $i=k+1,\cdots n$ we have $\lambda_i=\mu_i$).
So \be F(T_{z_2}R_{z_2^{-1}z_1}(Z_2))=F(Z_1).\ee But $F$ is a
right invariant Finsler metric on $G$, so for any $g_1, g_2\in G$
and $X_{g_1}\in T_{g_1}G$ we have \be
F(T_{g_1}R_{g_2}(X_{g_1}))=F(X_{g_1}),\nonumber \ee therefore \be
F(Z_2)=F(T_{z_2}R_{z_2^{-1}z_1}(Z_2)).\ee By equations 3.1 and
3.2 we have $F(Z_1)=F(Z_2)$.\\
It means the definition of ${\mathcal{F}}(X)$ is independent of
choice of $z$, so ${\mathcal{F}}$ is welldefined.\\
${\mathcal{F}}$ has all two conditions of Finsler metrics,
because ${\mathcal{F}}=F\circ (Tp|_V)^{-1}$. Also ${\mathcal{F}}$
is right invariant under right action of $G$ on $M$, because $F$
is right invariant on $G$.\hfill\ $\Box$
\paragraph{3.4. Remark.}
If we want to have a similar theorem as Theorem 3.3, for left
invariant Finsler metrics, it suffices to replace the word
``right'' by ``left'' in Lemma 3.2 and Theorem 3.3, and to use the
fact that $zH=Hz$ by the normality of $H$.

\paragraph{3.5. Theorem.}Assume that $G$ is any $n-$dimensional
Lie group, $H$ any closed normal Lie subgroup, $M=G/H$ the
quotient group and $p:G\rightarrow M$ is the natural projection.
If $F$ is a bi-invariant Finsler metric on $G$, then there is a
Finsler metric on $M$ induced by $F$ such that is invariant
under the natural right and left actions of $G$ on $M$.\\

\noindent \textit{Proof}: Suppose that $\frak{g}$, $\frak{h}$,
$\{X_1,\cdots,X_k,X_{k+1},\cdots,X_n\}$ and $V$ are the same
objects in the proof of Theorem 3.3. Let $z\in G$ and $Z\in V(z)$
such that $p(z)=x$ and $T_zp(Z)=X$. We define
${\mathcal{F}}:TM\rightarrow [0,\infty)$ by
${\mathcal{F}}(X):=F(Z)$. By Remark 3.4, this definition is
welldefined and also left invariant.\\
Since $\{X_1,\cdots,X_k\}$ is a basis of the Lie algebra consists
of left invariant vector fields of $H$, so $\{\nu^\ast
X_1,\cdots,\nu^\ast X_k\}$ is a basis of the Lie algebra consists
of right invariant vector fields of $H$ and $\{\nu^\ast
X_1,\cdots,\nu^\ast X_n\}$ is a basis of the Lie algebra consists
of right invariant vector fields of $G$. Now by using Lemma 3.2
and Theorem 3.3 and the fact that, $zH=Hz$ for any $z\in G$, we
have ${\mathcal{F}}$ is right invariant, therefore
${\mathcal{F}}$ is bi-invariant.\hfill\ $\Box$
\paragraph{3.6. Corollary.} Let $G$ be any $n-$dimensional
connected Lie group, $H$ any connected closed Lie subgroup and
$M:=G/H$ the quotient manifold. Suppose that $F$ is a left
invariant (right or bi-invariant) Finsler metric on $G$. If
$\frak{h}$ is an ideal of $\frak{g}$ then $M$ admits a left
invariant (right or bi-invariant) Finsler metric in the natural
way.\\

\noindent \textit{Proof}: Since $G$ and $H$ are connected Lie
groups and $\frak{h}$ is an ideal of $\frak{g}$, by  Theorem
2.13.4 of  \cite{[8]}, $H$ is a closed normal Lie subgroup of $G$.
So by attention to Theorems 3.3, 3.5 and Remark 3.4 the proof will
be finished.\hfill\ $\Box$
\paragraph{3.7. Corollary.} Let $G$ be any $n-$dimensional abelian
Lie group and $H$ a closed subgroup of $G$. If $F$ is a left
invariant (right or bi-invariant) Finsler metric on $G$ then $F$
induces a left invariant (right or bi-invariant) Finsler metric
on $M=G/H$ in the natural way.\\

Our results are true in Riemann case, and also our method for
construct invariant Finsler metrics on quotient groups is
compatible with method described in the first part of section 3
about invariant Riemannian metrics.

\end{document}